\documentstyle{amsppt}
\voffset-10mm
\magnification1200
\pagewidth{130mm}
\pageheight{204mm}
\hfuzz=2.5pt\rightskip=0pt plus1pt
\binoppenalty=10000\relpenalty=10000\relax
\TagsOnRight
\loadbold
\nologo
\addto\tenpoint{\normalbaselineskip=1.3\normalbaselineskip\normalbaselines}
\addto\eightpoint{\normalbaselineskip=1.2\normalbaselineskip\normalbaselines}

\let\ge\geqslant
\let\wt\widetilde
\define\pfrac#1#2{
\thickfracwithdelims..\thickness0{}{\lower2pt\rlap{\vrule height11.5pt}}
\kern-2pt
\frac{#1\,}{\,#2}
\kern-2.4pt
\thickfracwithdelims..\thickness0{\lower4.3pt\rlap{\vrule height11.5pt}}{}
}
\redefine\d{\roman d}

\redefine\Re{\operatorname{Re}}

\topmatter
\title
An Ap\'ery-like difference equation\\
for Catalan's constant
\endtitle
\author
Wadim Zudilin \rm(Moscow)
\endauthor
\date
\hbox to100mm{\vbox{\hsize=100mm%
\centerline{E-print \tt math.NT/0201024}
\smallskip
\centerline{18 January 2002}
}}
\enddate
\address
\hbox to70mm{\vbox{\hsize=70mm%
\leftline{Moscow Lomonosov State University}
\leftline{Department of Mechanics and Mathematics}
\leftline{Vorobiovy Gory, Moscow 119899 RUSSIA}
\leftline{{\it URL\/}: \tt http://wain.mi.ras.ru/index.html}
}}
\endaddress
\email
{\tt wadim\@ips.ras.ru}
\endemail
\abstract
Applying Zeilberger's algorithm of creative telescoping
to a family of certain very-well-poised hypergeometric
series involving linear forms in Catalan's constant
with rational coefficients, we obtain a second-order 
difference equation for these forms and their coefficients.
As a consequence we obtain a new way of fast calculation
of Catalan's constant as well as a new continued-fraction
expansion for it. Similar arguments can be put forward
to indicate a second-order difference equation and
a new continued fraction for $\zeta(4)=\pi^4/90$, and
we announce corresponding results at the end of this paper.
\endabstract
\keywords
Ap\'ery-like difference equation,
continued fraction, Catalan's constant, $\pi^4$,
generalized hypergeometric series,
Zeilberger's algorithm of creative telescoping
\endkeywords
\endtopmatter
\footnote""{2000 {\it Mathematics Subject Classification}.\enspace
Primary 39A05, 11B37, 11J70; Secondary 33C20, 33C60, 11B65, 11M06.}
\leftheadtext{W.~Zudilin}
\rightheadtext{Ap\'ery-like difference equation for Catalan's constant}
\document

\head
1. Introduction
\endhead
One of the most crucial and quite mysterious ingredients
in Ap\'ery's proof~\cite{Ap},~\cite{Po} of the irrationality
of~$\zeta(2)$ and~$\zeta(3)$ is the existence of the difference
equations
$$
\gathered
(n+1)^2u_{n+1}-(11n^2+11n+3)u_n-n^2u_{n-1}=0,
\\
u_0'=1, \quad u_1'=3, \qquad v_0'=0, \quad v_1'=5,
\endgathered
\tag1
$$
and
$$
\gathered
(n+1)^3u_{n+1}-(2n+1)(17n^2+17n+5)u_n+n^3u_n=0,
\\
u_0''=1, \quad u_1''=5, \qquad v_0''=0, \quad v_1''=6,
\endgathered
$$
with the following properties of their solutions:
$$
\lim_{n\to\infty}\frac{v_n'}{u_n'}=\zeta(2),
\qquad
\lim_{n\to\infty}\frac{v_n''}{u_n''}=\zeta(3).
$$
Unexpected inclusions
$u_n',D_n^2v_n'\in\Bbb Z$ and $u_n'',D_n^3v_n''\in\Bbb Z$,
where $D_n$~denotes
the least common multiple of the numbers $1,2,\dots,n$
(and $D_0=1$ for completeness), together with the prime number
theorem ($D_n^{1/n}\to e$ as $n\to\infty$) and Poincar\'e's theorem,
then yield the following asymptotic behaviour of the linear
forms $D_n^2u_n'\zeta(2)-D_n^2v_n'$ and
$D_n^3u_n''\zeta(3)-D_n^3v_n''$ with {\it integral\/} coefficients:
$$
\aligned
\lim_{n\to\infty}|D_n^2u_n'\zeta(2)-D_n^2v_n'|^{1/n}
&=e^2\biggl(\frac{\sqrt5-1}2\biggr)^5<1,
\\
\lim_{n\to\infty}|D_n^3u_n''\zeta(3)-D_n^3v_n''|^{1/n}
&=e^3(\sqrt2-1)^4<1,
\endaligned
$$
and thus one obtains that both $\zeta(2)$ and~$\zeta(3)$ cannot be rational.

The two following decades after \cite{Ap} were full of attempts
to indicate the total list of the second-order recursions with
integral solutions and to show their `geometric' (i.e., Picard--Fuchs
differential equations) origin~\cite{Za}. We do not pretend
to be so heroic in this paper, and we apply quite elementary
arguments to get new recurrence equations with `almost-integral'
solutions.

In our general, joint with T.~Rivoal, study~\cite{RZ}
of arithmetic properties for values of Dirichlet's beta function
$$
\beta(s):=\sum_{l=0}^{\infty}\frac{(-1)^l}{(2l+1)^s}
$$
at positive integers~$s$, we have discovered a construction
of $\Bbb Q$-linear forms in~$1$ and Catalan's constant
$$
G:=\sum_{l=0}^\infty\frac{(-1)^l}{(2l+1)^2}=\beta(2)
$$
similar to the one considered by Ap\'ery in his proof of
the irrationality of~$\zeta(2)$. The analogy is far from
proving the desired irrationality of~$G$, but it allows to indicate
the following second-order difference equation
$$
(2n+1)^2(2n+2)^2p(n)u_{n+1}-q(n)u_n-(2n-1)^2(2n)^2p(n+1)u_{n-1}=0,
\tag2
$$
where
$$
\aligned
p(n)
&=20n^2-8n+1,
\\
q(n)
&=3520n^6+5632n^5+2064n^4-384n^3-156n^2+16n+7,
\endaligned
\tag3
$$
with the initial data
$$
u_0=1, \quad u_1=\frac74, \qquad v_0=0, \quad v_1=\frac{13}8.
\tag4
$$

\proclaim{Theorem 1}
For each $n=0,1,2,\dots$, the numbers $u_n$ and $v_n$
produced by the recursion~\thetag{2},~\thetag{4}
are positive rationals satisfying the inclusions
$$
2^{4n+3}D_nu_n\in\Bbb Z, \qquad 2^{4n+3}D_{2n-1}^3v_n\in\Bbb Z,
\tag5
$$
and the following limit relation holds
$$
\lim_{n\to\infty}\frac{v_n}{u_n}=G.
$$
\endproclaim

The positivity and rationality of $u_n$ and $v_n$ follows
immediately from~\thetag{2}--\thetag{4}. The characteristic
polynomial $\lambda^2-11\lambda-1$ with
zeros $\bigl((1\pm\sqrt5\,)/2\bigr)^5$
of the difference equation~\thetag{2} coincides with
the corresponding one
for Ap\'ery's equation~\thetag{1}. Therefore
application of Poincar\'e's theorem
(see also \cite{Zu1, Proposition~2})
yields the limit relations
$$
\gather
\lim_{n\to\infty}u_n^{1/n}
=\lim_{n\to\infty}v_n^{1/n}
=\biggl(\frac{1+\sqrt5}2\biggr)^5=\exp(2.40605912\dots),
\\
\lim_{n\to\infty}|u_nG-v_n|^{1/n}
=\biggl|\frac{1-\sqrt5}2\biggr|^5=\exp(-2.40605912\dots),
\endgather
$$
while the inclusions~\thetag{5} and the prime number theorem
imply that the linear forms $2^{4n+3}D_{2n-1}^3(u_nG-v_n)$ with integral
coefficients do not tend to~$0$ as $n\to\infty$. However,
the rational approximations $v_n/u_n$ to Catalan's constant
converge quite rapidly (for instance, $|v_{10}/u_{10}-G|<10^{-20}$)
and one can use the recursion~\thetag{2},~\thetag{4} for
fast evaluating~$G$. Another consequence of Theorem~1 is a
new continued-fraction expansion for Catalan's constant.
Namely, considering $v_n/u_n$ as convergents of a continued
fraction for~$G$ and making the equivalent transform of the fraction
\cite{JT, Theorems 2.2 and 2.6} we arrive at

\proclaim{Theorem 2}
The following expansion holds:
$$
G=\pfrac{13/2}{q(0)}+\pfrac{1^4\cdot2^4\cdot p(0)p(2)}{q(1)}+\dots
+\pfrac{(2n-1)^4(2n)^4p(n-1)p(n+1)}{q(n)}+\dotsb,
$$
where the polynomials $p(n)$ and $q(n)$ are given in~\thetag{3}.
\endproclaim

The multiple-integral representation for the linear forms $u_nG-v_n$
similar to those obtained by F.~Beukers in \cite{Be, Formula~(5)}
for the linear forms $u_n'\zeta(2)-v_n'$ is given by

\proclaim{Theorem 3}
For each $n=0,1,2,\dots$, there holds the identity
$$
u_nG-v_n
=\frac{(-1)^n}4\int_0^1\!\!\int_0^1
\frac{x^{n-1/2}(1-x)^ny^n(1-y)^{n-1/2}}
{(1-xy)^{n+1}}\,\d x\,\d y.
\tag6
$$
\endproclaim

\head
2. Difference equation for Catalan's constant
\endhead

Consider the rational function
$$
R_n(t)
:=n!(2t+n+1)
\frac{t(t-1)\dotsb(t-n+1)\cdot(t+n+1)\dotsb(t+2n)}
{((t+\frac12)(t+\frac32)\dotsb(t+n+\frac12))^3}
\tag7
$$
and the corresponding (very-well-poised) hypergeometric series
$$
\align
F_n
&:=\sum_{t=0}^\infty(-1)^tR_n(t)
\\ &\phantom:
=(-1)^nn!\,
\frac{\Gamma(3n+2)\,\Gamma(n+\frac12)^3\Gamma(n+1)}
{\Gamma(2n+\frac32)^3\Gamma(2n+1)}
\\ &\phantom:\qquad\times
{}_6\!F_5\biggl(\matrix\format&\,\r\\
3n+1, & \frac32n+\frac32, &  n+\frac12, &  n+\frac12, &  n+\frac12, &  n+1 \\
      & \frac32n+\frac12, & 2n+\frac32, & 2n+\frac32, & 2n+\frac32, & 2n+1
\endmatrix\biggm|-1\biggr).
\tag8
\endalign
$$

\proclaim{Lemma 1}
There holds the equality
$$
F_n=U_n\beta(3)+U_n'\beta(2)+U_n''\beta(1)-V_n,
\tag9
$$
where $U_n,D_nU_n',D_n^2U_n'',D_{2n-1}^3V_n\in2^{-4n}\Bbb Z$.
\endproclaim

\demo{Proof}
We start with mentioning that
$$
P_n^{(1)}(t):=\frac{t(t-1)\dotsb(t-n+1)}{n!}
\qquad\text{and}\qquad
P_n^{(2)}(t):=\frac{(t+n+1)\dotsb(t+2n)}{n!}
\tag10
$$
are integral-valued polynomials and,
as it is known~\cite{Zu2, Lemma~7},
$$
2^{2n}\cdot P_n(-k-\tfrac12)\in\Bbb Z
\qquad\text{for}\quad k\in\Bbb Z
\tag11
$$
and, moreover,
$$
2^{2n}D_n^j\cdot\frac1{j!}\,\frac{\d^jP_n(t)}{\d t^j}\bigg|_{t=-k-1/2}
\in\Bbb Z
\qquad\text{for}\quad k\in\Bbb Z
\quad\text{and}\quad j=1,2,\dots,
\tag12
$$
where $P_n(t)$ is any of the polynomials~\thetag{10}.

The rational function
$$
Q_n(t):=\frac{n!}{(t+\frac12)(t+\frac32)\dotsb(t+n+\frac12)}
$$
has also `nice' arithmetic properties. Namely,
$$
a_k:=Q_n(t)(t+k+\tfrac12)\big|_{t=-k-1/2}=\cases
(-1)^k\binom nk\in\Bbb Z & \text{if $k=0,1,\dots,n$}, \\
0 & \text{for other $k\in\Bbb Z$},
\endcases
\tag13
$$
that allow to write the following partial-fraction expansion:
$$
Q_n(t)=\sum_{l=0}^n\frac{a_l}{t+l+\frac12}.
$$
Hence, for $j=1,2,\dots$ we obtain
$$
\align
\frac{D_n^j}{j!}\,\frac{\d^j}{\d t^j}
\bigl(Q_n(t)(t+k+\tfrac12)\bigr)\big|_{t=-k-1/2}
&=\frac{D_n^j}{j!}\,\frac{\d^j}{\d t^j}
\sum_{l=0}^na_l\biggl(1-\frac{l-k}{t+l+\frac12}\biggr)\bigg|_{t=-k-1/2}
\\
&=(-1)^{j-1}D_n^j\sum\Sb l=0\\l\ne k\endSb^n\frac1{(l-k)^j}\in\Bbb Z.
\tag14
\endalign
$$
Therefore the inclusions \thetag{11}--\thetag{14}
and the Leibniz rule for differentiating a product
imply that the numbers
$$
\align
A_{jk}
&=A_{jk}(n)
:=\frac1{j!}\,\frac{\d^j}{\d t^j}
\bigl(R_n(t)(t+k+\tfrac12)^3\bigr)\big|_{t=-k-1/2}
\tag15
\\
&=\frac1{j!}\,\frac{\d^j}{\d t^j}
\bigl((2t+n+1)\cdot P_n^{(1)}(t)\cdot P_n^{(2)}(t)
\cdot(Q_n(t)(t+k+\tfrac12))^3\bigr)\big|_{t=-k-1/2}
\endalign
$$
satisfy the inclusions
$$
2^{4n}D_n^j\cdot A_{jk}(n)\in\Bbb Z
\qquad\text{for}\quad k=0,1,\dots,n
\quad\text{and}\quad j=0,1,2,\dots\,.
\tag16
$$
Mention now that the numbers \thetag{15} are coefficients
in the partial-fraction expansion of the rational function~\thetag{7},
$$
R_n(t)=\sum_{j=0}^2\sum_{k=0}^n\frac{A_{jk}}{(t+k+\frac12)^{3-j}}.
\tag17
$$
Substituting this expansion into the definition~\thetag{8}
of the quantity~$F_n$ we obtain the desired representaion~\thetag{9}:
$$
\align
F_n
&=\sum_{j=0}^2\sum_{k=0}^n(-1)^kA_{jk}
\sum_{t=0}^\infty\frac{(-1)^{t+k}}{(t+k+\frac12)^{3-j}}
\\
&=\sum_{j=0}^2\sum_{k=0}^n(-1)^kA_{jk}
\biggl(\sum_{l=0}^\infty-\sum_{l=0}^{k-1}\biggr)
\frac{(-1)^l}{(l+\frac12)^{3-j}}
\\
&=U_n\beta(3)+U_n'\beta(2)+U_n''\beta(1)-V_n,
\endalign
$$
where
$$
\gather
U_n=2^3\sum_{k=0}^n(-1)^kA_{0k}(n),
\quad
U_n'=2^2\sum_{k=0}^n(-1)^kA_{1k}(n),
\quad
U_n''=2\sum_{k=0}^n(-1)^kA_{2k}(n),
\tag18
\\
V_n=\sum_{j=0}^22^{3-j}\sum_{k=0}^n(-1)^kA_{jk}(n)
\sum_{l=0}^{k-1}\frac{(-1)^l}{(2l+1)^{3-j}}.
\tag19
\endgather
$$
Finally, using the inclusions~\thetag{16} and
$$
D_{2n-1}^{3-j}\sum_{l=0}^{k-1}\frac{(-1)^l}{(2l+1)^{3-j}}\in\Bbb Z
\qquad\text{for}\quad k=0,1,\dots,n
\quad\text{and}\quad j=0,1,2,
$$
we deduce that $U_n,D_nU_n',D_n^2U_n'',D_{2n-1}^3V_n\in2^{-4n}\Bbb Z$
as required.
\enddemo

Using Zeilberger's algorithm of creative telescoping
\cite{PWZ, Section~6} for the rational function~\thetag{7},
we obtain the certificate $S_n(t):=s_n(t)R_n(t)$, where
$$
\align
&
s_n(t)
:=\frac1{2(2t+n+1)(t+2n-1)(t+2n)}\cdot
\bigl(8n(2n-1)^2(20n^2+32n+13)t^4
\\ &\quad
+2(5440n^6+7104n^5+912n^4-1088n^3+76n^2+68n+7)t^3
\\ &\quad
+(44800n^7+65600n^6+17568n^5-7056n^4-1088n^3+372n^2+146n-1)t^2
\\ &\quad
+(2n+1)(34880n^7+39328n^6-2176n^5-8416n^4+964n^3+154n^2+58n-13)t
\\ &\quad
+n(2n-1)(2n+1)^2(4720n^5+6192n^4+816n^3-864n^2+69n+13)\bigr)
\tag20
\endalign
$$
satisfying the following property.

\proclaim{Lemma 2}
For each $n=1,2,\dots$, there holds the identity
$$
\align
&
(2n+1)^2(2n+2)^2p(n)R_{n+1}(t)-q(n)R_n(t)
-(2n-1)^2(2n)^2p(n+1)R_{n-1}(t)
\\ &\qquad
=-S_n(t+1)-S_n(t),
\tag21
\endalign
$$
where the polynomials $p(n)$ and $q(n)$ are given in~\thetag{3}.
\endproclaim

\demo{Proof}
Divide both sides of~\thetag{21} by $R_n(t)$ and verify
the identity
$$
\align
&
(2n+1)^2(2n+2)^2p(n)\cdot(n+1)
\frac{(2t+n+2)(t-n)(t+2n+1)(t+2n+2)}{(2t+n+1)(t+n+1)(t+n+\frac32)^3}
-q(n)
\\ &\quad\qquad
-(2n-1)^2(2n)^2p(n+1)\cdot
\frac{(2t+n)(t+n)(t+n+\frac12)^3}{n(2t+n+1)(t-n+1)(t+2n-1)(t+2n)}
\\ &\quad
=-s_n(t+1)\frac{(2t+n+3)(t+\frac12)^3(t+1)(t+2n+1)}
{(2t+n+1)(t-n+1)(t+n+1)(t+n+\frac32)^3}
-s_n(t),
\endalign
$$
where $s_n(t)$ is given in~\thetag{20}.
\enddemo

\proclaim{Lemma 3}
The quantity~\thetag{8} satisfies the difference equation~\thetag{2}
for $n=1,2,\dots$\,.
\endproclaim

\demo{Proof}
Multiplying both sides of the equality~\thetag{21} by~$(-1)^t$
and summing the result over $t=0,1,2,\dots$
we obtain
$$
(2n+1)^2(2n+2)^2p(n)F_{n+1}-b(n)F_n-(2n-1)^2(2n)^2p(n+1)F_{n-1}
=-S_n(0).
$$
It remains to note that, for $n\ge1$, both functions
$R_n(t)$ and $S_n(t)=s_n(t)R_n(t)$ have zero
at $t=0$. Thus $S_n(0)=0$ for $n=1,2,\dots$ and we obtain
the desired recurrence~\thetag{2} for the quantity~\thetag{8}.
\enddemo

\proclaim{Lemma 4}
The coefficients $U_n,U_n',U_n'',V_n$ in the representation
\thetag{9} satisfy the difference equation~\thetag{2}
for $n=1,2,\dots$\,.
\endproclaim

\demo{Proof}
We can write down the partial-fraction expansion~\thetag{17} in the form
$$
R_n(t)=\sum_{j=1}^4\sum_{k=-\infty}^{+\infty}
\frac{A_{jk}(n)}{(t+k+\frac12)^{3-j}},
$$
where the formulae~\thetag{15} remain true for {\it all\/} $k\in\Bbb Z$
(not for $k=0,1,\dots,n$ only). Now,
multiply both sides of~\thetag{21} by $(-1)^k(t+k+\frac12)^3$,
take the $j$th derivative, where $j=0,1,2$,
substitute $t=-k-\frac12$ in the result,
and sum over all integers~$k$;
this procedure implies that the numbers
$$
U_n=8\sum_{k=-\infty}^{+\infty}(-1)^kA_{0k}(n),
\quad
U_n'=4\sum_{k=-\infty}^{+\infty}(-1)^kA_{1k}(n),
\quad
U_n''=2\sum_{k=-\infty}^{+\infty}(-1)^kA_{2k}(n)
$$
(cf.~\thetag{18}) satisfy the difference equation~\thetag{2}.
Finally, the sequence
$$
V_n=U_n\beta(3)+U_n'\beta(2)+U_n''\beta(1)-F_n
$$
also satisfies the recursion~\thetag{2} by Lemma~3 and the above.
\enddemo

Since
$$
R_0(t)=\frac2{(t+\frac12)^2},
\qquad
R_1(t)
=-\frac{3/4}{(t+\frac12)^3}-\frac{3/4}{(t+\frac32)^3}
+\frac{7/4}{(t+\frac12)^2}-\frac{7/4}{(t+\frac32)^2},
$$
in accordance with~\thetag{18},~\thetag{19} we obtain
$$
U_0'=8, \quad U_0=U_0''=V_0=0,
\qquad\text{and}\qquad
U_1'=14, \quad V_1=13, \quad U_1=U_1''=0,
$$
hence as a consequence of Lemma~4 we deduce that $U_n=U_n''=0$
for $n=0,1,2,\dots$\,.

\proclaim{Lemma 5}
There holds the equality
$$
F_n=U_n'G-V_n,
$$
where $2^{4n}D_nU_n'\in\Bbb Z$ and $2^{4n}D_{2n-1}^3V_n\in\Bbb Z$.
\endproclaim

The sequences $u_n:=U_n'/8$ and $v_n:=V_n/8$
satisfy the difference equation~\thetag{2}
and initial conditions~\thetag{4};
the fact that $F_n\ne0$ and $|F_n|\to0$ as $n\to\infty$
follows from Theorem~3 and asymptotics of the multiple integral~\thetag{6}
proved in~\cite{Be}. This completes the proof of Theorem~1.

\head
3. Connection with ${}_3F_2$-hypergeometric series
\endhead

The corresponding very-well-poised hypergeometric series~\thetag{8}
at $z=-1$ can be reduced to a simpler series with the help
of Bailey's identity \cite{Ba, Section~4.4, formula~(2)}:
$$
\multline
{}_3\!F_2\biggl(\gathered
1+a-b-c, \, d, \, e \\ 1+a-b, \, 1+a-c
\endgathered\biggm|1\biggr)
=\frac{\Gamma(1+a)\,\Gamma(1+a-d-e)}{\Gamma(1+a-d)\,\Gamma(1+a-e)}
\\
\times{}_6\!F_5\biggl(\matrix\format&\,\r\\
a, & 1+\frac12a, &     b, &     c, &     d, &     e \\
   &   \frac12a, & 1+a-b, & 1+a-c, & 1+a-d, & 1+a-e
\endmatrix\biggm|-1\biggr),
\endmultline
$$
if $\Re(1+a-d-e)>0$. Namely,
in the case $a=3n+1$, $b=c=d=n+\frac12$, and $e=n+1$,
we obtain
$$
F_n=U_n'G-V_n
=(-1)^n\cdot2\int_0^1\!\!\int_0^1
\frac{x^{n-1/2}(1-x)^ny^n(1-y)^{n-1/2}}
{(1-xy)^{n+1}}\,\d x\,\d y,
$$
where the Euler-type integral representation for
the ${}_3\!F_2$-series can be derived
as in \cite{Sl, Section~4.1} and \cite{BR, proof of Lemma~2}.

This completes the proof of Theorem~3.

\head
4. Conclusion remarks
\endhead

The conclusion~\thetag{5} of Theorem~1 is far from being precise;
in fact, using \thetag{2},~\thetag{4} one gets experimentally
(up to $n=1000$, say) the stronger inclusions\footnote{%
A slightly weakened form of the inclusions is proved in~\cite{Zu4}.}
$$
2^{4n}u_n\in\Bbb Z, \qquad 2^{4n}D_{2n-1}^2v_n\in\Bbb Z.
$$
Unfortunately, they also give no chance to prove that Catalan's
constant is irrational since linear forms
$2^{4n}D_{2n-1}^2(u_nG-v_n)$ do not tend to~$0$ as $n\to\infty$.

In the same vein, using another very-well-poised series
of hypergeometric type
$$
\align
\wt F_n
&:=\frac{(-1)^{n+1}}6\sum_{t=1}^\infty
\frac{\d}{\d t}\biggl((2t+n)\frac{((t-1)\dotsb(t-n))^2
\cdot((t+n+1)\dotsb(t+2n))^2}
{(t(t+1)\dotsb(t+n))^4}\biggr)
\\ &\phantom:
=\wt u_n\zeta(4)-\wt v_n
\endalign
$$
and the arguments of Section~2, we deduce the difference equation
$$
(n+1)^5u_{n+1}-r(n)u_n-3n^3(3n-1)(3n+1)u_{n-1}=0,
\tag22
$$
where
$$
\align
r(n)
&=3(2n+1)(3n^2+3n+1)(15n^2+15n+4)
\\
&=270n^5+675n^4+702n^3+378n^2+105n+12,
\tag23
\endalign
$$
with the initial data
$$
\wt u_0=1, \quad \wt u_1=12, \qquad \wt v_0=0, \quad \wt v_1=13
$$
for its two independent solutions $\wt u_n$ and $\wt v_n$,

\proclaim{Theorem 4 \cite{Zu3}}
For each $n=0,1,2,\dots$, the numbers $\wt u_n$ and $\wt v_n$
are positive rationals satisfying the inclusions
$$
6D_n\wt u_n\in\Bbb Z, \qquad 6D_n^5\wt v_n\in\Bbb Z,
\tag24
$$
and there holds the limit relation
$$
\lim_{n\to\infty}\frac{\wt v_n}{\wt u_n}=\frac{\pi^4}{90}
=\zeta(4)=\sum_{n=1}^\infty\frac1{n^4}.
\tag25
$$
\endproclaim

\remark{Remark}
During the preparation of this article,
we have known that the difference equation~\thetag{22},
in slightly different normalization,
and the limit relation~\thetag{25} without the inclusions~\thetag{24}
had been stated independently
by H.~Cohen and G.~Rhin~\cite{Co} using Ap\'ery's
`acc\'el\'eration de la convergence' approach,
and by V.\,N.~Sorokin~\cite{So} by means
of certain explicit Hermite--Pad\'e-type approximations.
We underline that our approach differs from that
of~\cite{Co} and~\cite{So}.
\endremark

Application of Poincar\'e's theorem yields the asymptotic relations
$$
\lim_{n\to\infty}|\wt u_n|^{1/n}
=\lim_{n\to\infty}|\wt v_n|^{1/n}
=(3+2\sqrt3\,)^3=\exp(5.59879212\dots)
$$
and
$$
\lim_{n\to\infty}|\wt u_n\zeta(4)-\wt v_n|^{1/n}
=|3-2\sqrt3\,|^3=\exp(-2.30295525\dots),
$$
since the characteristic polynomial $\lambda^2-270\lambda-27$
of the equation~\thetag{22}
has zeros $135\pm78\sqrt3=(3\pm2\sqrt3\,)^3$.
Thus, we can consider $\wt v_n/\wt u_n$ as convergents of
a continued fraction for~$\zeta(4)$ and making
the equivalent transform of the fraction we obtain

\proclaim{Theorem 5}
The following expansion holds:
$$
\zeta(4)=\pfrac{13}{r(0)}+\pfrac{1^7\cdot2\cdot3\cdot4}{r(1)}
+\pfrac{2^7\cdot5\cdot6\cdot7}{r(2)}+\dots
+\pfrac{n^7(3n-1)(3n)(3n+1)}{r(n)}+\dotsb,
$$
where the polynomial $r(n)$ is given in~\thetag{23}.
\endproclaim

The mystery of the $\zeta(4)$-case consists in the fact
that experimental calculations give us the better inclusions
$$
\wt u_n\in\Bbb Z, \qquad D_n^4\wt v_n\in\Bbb Z
$$
(cf.~\thetag{24}); unfortunately, the linear forms
$D_n^4(\wt u_n\zeta(4)-\wt v_n)$
do not tend to~$0$ as $n\to\infty$.



\Refs
\widestnumber\key{PWZ}

\ref\key Ap
\by R.~Ap\'ery
\paper Irrationalit\'e de $\zeta(2)$ et $\zeta(3)$
\jour Ast\'erisque
\vol61
\yr1979
\pages11--13
\endref

\ref\key Ba
\by W.\,N.~Bailey
\book Generalized hypergeometric series
\bookinfo Cambridge Math. Tracts
\vol32
\publ Cambridge Univ. Press
\publaddr Cambridge
\yr1935
\moreref
\bookinfo 2nd reprinted edition
\publaddr New York--London
\publ Stechert-Hafner
\yr1964
\endref

\ref\key BR
\by K.~Ball and T.~Rivoal
\paper Irrationalit\'e d'une infinit\'e de valeurs
de la fonction z\^eta aux entiers impairs
\jour Invent. Math.
\yr2001
\vol146
\issue1
\pages193--207
\endref

\ref\key Be
\by F.~Beukers
\paper A note on the irrationality of~$\zeta(2)$ and~$\zeta(3)$
\jour Bull. London Math. Soc.
\vol11
\issue3
\yr1979
\pages268--272
\endref

\ref\key Co
\by H.~Cohen
\paper Acc\'el\'eration de la convergence
de certaines r\'ecurrences lin\'eaires
\inbook S\'eminaire de Th\'eorie des Nombres de Grenoble
(6 et 13 novembre 1980), 47 pages
\endref

\ref\key JT
\by W.\,B.~Jones and W.\,J.~Thron
\book Continued fractions. Analytic theory and applications
\bookinfo Encyclopaedia Math. Appl. Section: Analysis
\vol11
\publ Addison-Wesley
\publaddr London
\yr1980
\endref

\ref\key PWZ
\by M.~Petkov\v sek, H.\,S.~Wilf, and D.~Zeilberger
\book $A=B$
\publaddr Wellesley, M.A.
\publ A.\,K.~Peters, Ltd.
\yr1997
\endref

\ref\key Po
\by A.~van der Poorten
\paper A proof that Euler missed...
Ap\'ery's proof of the irrationality of~$\zeta(3)$
\paperinfo An informal report
\jour Math. Intelligencer
\vol1
\issue4
\yr1978/79
\pages195--203
\endref

\ref\key RZ
\by T.~Rivoal and W.~Zudilin
\paper Diophantine properties of numbers related to Catalan's constant
\inbook Pr\'epubl. de l'Institut de Math. de Jussieu,
no.~315 (January 2002), 17~pages
\endref

\ref\key Sl
\by L.\,J.~Slater
\book Generalized hypergeometric functions
\bookinfo 2nd edition
\publ Cambridge Univ. Press
\publaddr Cambridge
\yr1966
\endref

\ref\key So
\by V.\,N.~Sorokin
\paper One algorithm for fast calculation of~$\pi^4$
\inbook Preprint (April 2002)
\publaddr Moscow
\publ Russian Academy of Sciences,
M.\,V.~Kel\-dysh Institute for Applied Mathematics 
\yr2002
\endref

\ref\key Za
\by D.~Zagier
\book Integral solutions of Ap\'ery-like recurrence equations
\bookinfo Preprint
\yr2001
\endref

\ref\key Zu1
\by W.~Zudilin
\paper Difference equations and the irrationality measure of numbers
\paperinfo Analytic number theory and applications, Collection of papers
\jour Trudy Mat\. Inst\. Steklov [Proc\. Steklov Inst\. Math.]
\vol218
\yr1997
\pages165--178
\endref

\ref\key Zu2
\by W.~Zudilin
\paper Cancellation of factorials
\jour Mat. Sb. [Russian Acad. Sci. Sb. Math.]
\vol192
\yr2001
\issue8
\pages95--122
\moreref
\inbook E-print {\tt math.NT/0008017}
\endref

\ref\key Zu3
\by W.~Zudilin
\paper Well-poised hypergeometric service
for diophantine problems of zeta values
\paperinfo Actes des 12\`emes rencontres arithm\'e\-tiques de Caen
(June 29--30, 2001)
\jour J. Th\'eorie Nombres Bordeaux
\yr2003
\toappear
\endref

\ref\key Zu4
\by W.~Zudilin
\paper A few remarks on linear forms involving Catalan's constant
\inbook E-print {\tt math.NT/\allowlinebreak0210423} (21 October 2002)
\endref

\endRefs
\enddocument